\documentclass{conm-p-l}
\usepackage{amsmath}
\usepackage{amsthm}
\usepackage{amsfonts}
\usepackage{amssymb}
\usepackage{eucal}
\usepackage{graphicx}
\usepackage[all,knot]{xy}
\xyoption{arc}

\newcommand{\CC}{\mathcal{C}}
\newcommand{\sZ}{\mathcal{Z}}
\newcommand{\sS}{\mathcal{S}}
\newcommand{\fsl}{\mathfrak{sl}}

\newcommand{\End}{{\rm End}}
\newcommand{\Rep}{{\rm Rep}}
\newcommand{\U}{{\rm U}}
\newcommand{\SU}{{\rm SU}}

\newcommand{\ot}{\otimes}
\newcommand{\B}{\mathcal{B}}
\newcommand{\PP}{\mathcal{P}}

\newcommand{\lan}{\langle}
\newcommand{\ra}{\rangle}
\newcommand{\Hom}{{\rm Hom}}
\newcommand{\ve}{\varepsilon}

\newcommand{\Z}{\mathbb{Z}}
\newcommand{\g}{\mathfrak{g}}
\newcommand{\N}{\mathbb{N}}

\newcommand{\C}{\mathbb{C}}

\newcommand{\al}{\alpha}
\newtheorem{theorem}{Theorem}[section]

\theoremstyle{definition}
\newtheorem{definition}[theorem]{Definition}
\newtheorem{example}[theorem]{Example}
\newtheorem{question}[theorem]{Question}
\newtheorem{conjecture}[theorem]{Conjecture}
\theoremstyle{remark}
\newtheorem{remark}[theorem]{Remark}

\numberwithin{equation}{section}

\begin{document}

\title{Two paradigms for topological quantum computation}

\author{Eric C. Rowell}
\address{Department of Mathematics, Texas A\& M University
College Station, TX 77845}

\email{rowell@math.tamu.edu}

\subjclass[2000]{Primary 81P68; Secondary 20F36, 57M25, 68Q17, 57M27}
\date{March 8, 2008.}

\keywords{quantum computation, braid group, modular category, link invariant}
\begin{abstract}We present two paradigms relating algebraic, topological and quantum computational statistics for the topological
model for quantum computation.  In particular we suggest
correspondences between the computational power of topological quantum
computers, computational complexity of link invariants and images of braid group representations.
 While at least parts of these paradigms are well-known to experts, we provide supporting evidence for them in terms of recent results.  We give a fairly comprehensive list of known examples and formulate two
conjectures that would further support the paradigms.
\end{abstract}

\maketitle

\section{Introduction}

\begin{figure}[t0]

\centerline{\includegraphics[width=3.45in]{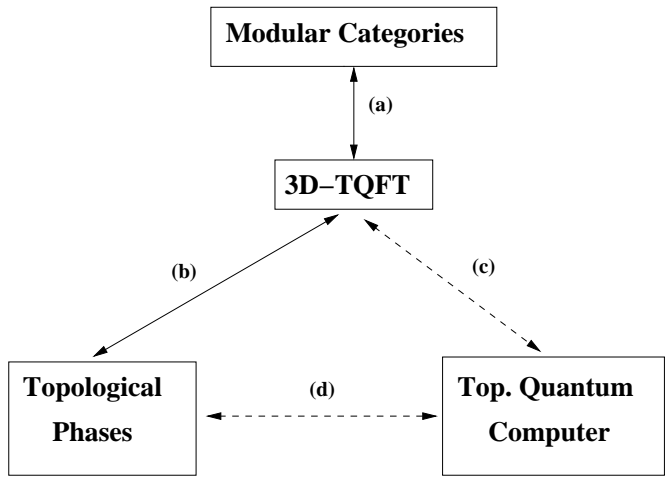}}

\caption{(a) equivalent by \cite{Tur}, (b) essentially by definition, see \cite{dSetal}, (c) idea
originated in \cite{FPNAS}, (d) first described in \cite{Kit} }
\label{fig:1}
\end{figure}

Topological quantum computation (TQC) is expected to be physically realized on quantum systems in \emph{topological phases}.  For example, the quasi-particle excitations in fractional quantum hall liquids are conjectured to exhibit the topological behavior necessary to support TQC.  A definition of topological phase is found in \cite{dSetal}: ``\textit{...a system is in a topological phase if its low-energy effective field theory is a topological quantum field theory (TQFT)}''.
Thus all observable properties of topological phases should be expressible in terms of the structure of the corresponding TQFT.  On the other hand, it is known \cite{Tur} that \emph{modular categories} faithfully encode (3D) TQFTs in algebraic terms.  These relationships between modular categories, TQFT, topological phases
and topological quantum computers are illustrated in Figure \ref{fig:1}.  The solid arrows represent
well-established or (tautological) one-to-one correspondences, while the dashed arrows represent theoretical expectations.

While the algebraic axioms defining modular categories may
seem quite distant from condensed matter physics and quantum computation, certain natural statistics in modular categories appear to correspond to important computational properties in TQC.  Most significantly, the images of the
braid group representations associated to a modular category are intimately
related to the computational power of the corresponding TQC.  We illustrate this with two well-known examples:
\begin{example}
 Consider the (unitary) modular category $\CC(\fsl_2,e^{\pi i/5})$ obtained as
a subquotient of the representation category of the quantum group $U_q\fsl_2$ with $q=e^{\pi i/5}$ (associated
with the $\SU(2)$-Chern-Simons-Witten TQFT at level $3$).  We note the following:
\begin{enumerate}
 \item The images of the associated braid group representations are as large as
possible, i.e. \emph{dense} in the group of special unitaries \cite{FLW2}.
\item A topological quantum computer realized upon a physical system algebraically modeled by $\CC(\fsl_2,e^{\pi i/5})$ is \textit{universal}.
\item The associated link invariant is $J_L(e^{2\pi i/5})$ the Jones polynomial evaluated at $q^2$, which has computational complexity $\#P$-hard \cite{JVW}.

\end{enumerate}
\end{example}
Moreover, approximate computation of the Jones polynomial at a $5$th root of unity is known to be
$BQP$-complete, so that it is essentially the hardest problem any quantum computer can hope to solve.
\begin{example}
Consider the (unitary) modular category $\CC(\fsl_2,e^{\pi i/4})$ obtained as
a subquotient of the representation category of the quantum group $U_q\fsl_2$ with $q=e^{\pi i/4}$ (associated
with the $\SU(2)$-Chern-Simons-Witten TQFT at level $2$).  We note the following:
\begin{enumerate}
 \item The images of the associated braid group representations factor over
finite groups \cite{jones86}.
\item A topological quantum computer realized upon a physical system algebraically modeled by $\CC(\fsl_2,e^{\pi i/4})$, while highly-entangling (see \cite{FRW}, \cite{KL}) is not universal.
\item The associated link invariant is $J_L(i)$ the Jones polynomial evaluated at $q^2=i$, which can be computed in polynomial time \cite{JVW}.

\end{enumerate}
\end{example}

Two paradigms based upon these examples might then associate density of
the braid group image with universal quantum computers and $\#P$-hard computational problems, and finite braid group images with non-universal (but potentially entangling) quantum devices and polynomial-time
computational problems.  Indeed, all evaluations of the Jones, HOMFLYPT and Kauffman  link invariants that are polynomial-time computable on a classical computer are associated with ``classical'' link invariants (see \cite[Theorems 6.3.2, 6.3.5 and 6.3.6]{Welsh}), by which we
mean link invariants pre-dating quantum topology.  Moreover, the corresponding braid group images in these cases have been shown to be finite in essentially all cases (see \cite{FLW2}, \cite{GoldJones}, \cite{LRW}, \cite{LR}, \cite{jones86}).
However, a deeper examination of further examples reveals
that this is not quite correct and a slight refinement is necessary.

In this paper we intend to describe such a refinement of these two paradigms.
Our aim is three-fold: to give paradigms that
can be of theoretical value to physicists, to describe a few conjectures of mathematical interest, and to present one perspective on the landscape of inter-related fields represented in the topological quantum computation endeavor.

Any attempt to be fully self-contained would require the introduction of many concepts from category theory, low-dimensional topology, complexity theory and condensed matter physics.  For brevity's sake, we
will content ourselves with providing the reader with a few references.  For an excellent survey of the physical and theoretical
set-up for TQCs, see \cite{dSetal}.  For the categorical and topological concepts, see \cite{BK} and \cite{Tur}.  For complexity theory applied to topological invariants, see \cite{Welsh}.

\subsubsection*{Acknowledgments}
The author would like to thank the following people for their generosity in valuable correspondence
and conversations: S. Witherspoon, L. Goldberg, Y. Zhang, J. Ospina,
M. Rojas, A. Bulatov, Z. Wang, G. Kuperberg, T. Stanford, and M. Thistlethwaite.

\section{Background}
We briefly describe some of the important features of modular categories and their relationships with topological phases, link invariants and
topological quantum computers.

A unitary modular category (UMC) $\CC$ is a semisimple $\C$-linear rigid ribbon category of finite rank satisfying
a certain non-degeneracy condition, such that the morphism spaces are equipped with a positive definite
hermitian form compatible with the other structures.  The representation category of a finite group is an
example of a category that satisfies all but one of the defining axioms of unitary modular tensor categories: namely it fails the \emph{modularity} (non-degeneracy) condition.  UMCs are constructed in a diversity of ways from
various fields of mathematics.  

\subsection{Constructions of UMCs} Often very different constructions yield equivalent categories, so we will
only list a few well-known explicit constructions.
\begin{enumerate}
\item \textbf{Quantum groups.}
To any finite dimensional simple Lie algebra $\g$ and a root of unity $q=e^{\pi i/\ell}$
one may associate a pre-modular category $\CC(\g,q)$.  These are obtained as subquotients of
the category of finite dimensional representations of the quantum group $U_q\g$, see \cite{Rsurvey} for
a survey.  Such a category may fail to be modular
or unitary (see \cite{Rowell1} and \cite{RJPAA}), but such circumstances can be avoided by certain restrictions on $\ell$.  Specifically, define $m=1$ for Lie types $A,D$ and $E$, $m=2$ for Lie types $B,C$ and $F_4$ and $m=3$ for Lie type $G_2$.  Then $\CC(\g,q)$ is a UMC provided $m\mid\ell$ (see \cite{Wenzlcstar}).
\item \textbf{Finite groups.}  Fix a finite group $G$ and a 3-cocycle $\omega$.  Then the twisted
 double of $G$, $D^\omega G$, is a finite dimensional quasi-triangular quasi-Hopf algebra.  The
representation category $\Rep(D^\omega G)$ is always a UMC, see \cite{BK} for details.
\item \textbf{Doubled spherical categories.}
There is a doubling procedure from which one obtains a modular category $\sZ(\sS)$ from a \emph{spherical category} $\sS$ (see \cite{BarWest} for the precise definition, and \cite{Muger2} for the double construction).  Briefly, a spherical category is a tensor category that is not necessarily braided but for which one has a canonical trace function.  Examples are ribbon categories and certain categories obtained from von Neumann algebras (see \cite{Izumi} for a description of the latter).  In fact, the representation categories of twisted doubles of finite groups can be obtained as the double of the spherical category $\Rep(\C[G])$ of representations of the group algebra of $G$.  Very few explicit ``new'' examples of modular categories obtained in this way have been worked out.  A few infinite families can be found in \cite{Izumi}, and the analysis of two examples are worked out in detail in \cite{HRW}.  If the spherical category $\sS$ is unitary the double $\sZ(\sS)$ will be a UMC.
\end{enumerate}

\subsection{Braid Group Representations}

The axioms of a UMC imply that for any object $X$ in a UMC $\CC$ one obtains a (highly non-degenerate) unitary
representation $\phi_X^n:\B_n\rightarrow\U(\End(X^{\ot n}))$.  Recall that $\B_n$, the braid group on $n$-strands, is the group with $n-1$ generators $\sigma_1,\ldots,\sigma_{n-1}$ satisfying:
\begin{enumerate}
 \item[(B1)] $\sigma_i\sigma_j=\sigma_j\sigma_i$ if $|i-j|\ge 2$
\item[(B2)] $\sigma_i\sigma_{i+1}\sigma_i=\sigma_{i+1}\sigma_i\sigma_{i+1}$ for $1\leq i\leq n-2$.
\end{enumerate}
The braiding on $\CC$ requires
that there is are natural braiding isomorphisms $C_{X,Y}:X\ot Y\cong Y\ot X$.  In particular
one obtains natural isomorphisms
$$R_X^i:=Id_X^{\ot (i-1)}\ot C_{X,X}\ot Id_X^{\ot (n-i-1)}\in\End(X^{\ot n})$$
so that the left action of $\End(X^{\ot n})$ on itself induces the representation $\phi_X^n$
by $$\phi_X^n(\sigma_i)f=R_X^i\circ f.$$
The unitarity of $\phi_X^n$ is due to the fact that $\End(X^{\ot n})$ is a Hilbert space, the naturality
of the braiding isomorphisms and the compatibility of the hermitian form with the other structures.

Given such a representation it is natural to ask
\begin{question}
 What is the closure of $\phi_X^n(\B_n)$ in $\U(\End(X^{\ot n}))$?
\end{question}
Indeed, this question was asked by Jones in \cite{jones86} long before its relevance to quantum computing
was realized.

Let us suppose that we have a decomposition $\End(X^{\ot n})=\bigoplus_k V_k$ into irreducible $\B_n$-representations,
and fix one irreducible subrepresentation $V_k$.
Denote by
$\Gamma_k$ the closure of the image of $\B_n$ in $\U(V)$.
Then $\Gamma_k$ modulo its center is exactly one of the following:
\begin{enumerate}
\item A finite abelian group
 \item A finite non-abelian group
\item An infinite compact group containing $\SU(V)$
\item An infinite compact group not containing $\SU(V)$
\end{enumerate}

These motivate the following:
\begin{definition}
\begin{enumerate}
 \item
 If $\Gamma_k/Z(\Gamma_k)$ is always a finite group for all objects $X$ in $\CC$, all
$n\in\N$, and all irreducible subrepresentations $V_k\subset\End(X^{\ot n})$ then we say $\CC$ has \textbf{property F}.
\item If $\Gamma_k/Z(\Gamma_k)$ is always a finite \emph{abelian} group for all objects $X$ in $\CC$, all
$n\in\N$, and all irreducible subrepresentations $V_k\subset\End(X^{\ot n})$ then we say $\CC$ has \textbf{property A}.  (Observe that this is the case whenever $\dim V_k=1$ for all $X,n$ and $k$.)
\item
 If there exists an object $X$ in $\CC$ and $N\in\N$ such that for all $n\geq N$ and for each irreducible subrepresentation $V_k\subset\End(X^{\ot n})$ the group $\Gamma_k/Z(\Gamma_k)$ contains $\SU(V_k)$ we say $\CC$ has the \textbf{density property}.
\end{enumerate}
\end{definition}

In nearly all cases one encounters in the literature, $\CC$ has either property \textbf{F} or the density property (see e.g. \cite{jones86}, \cite{FLW2}, \cite{LRW}, \cite{ERW} and \cite{LR}).
\begin{remark}
 One may generalize the construction above in the following way.  The \emph{pure braid group} $\PP_n$ is the (normal) subgroup of $\B_n$ generated by the conjugacy class of $\sigma_1^2$, or equivalently, the kernel of the obvious homomorphism $\B_n\rightarrow S_n$ that sends $\sigma_i$ to the transposition
$(i,i+1)$.  So geometrically $\PP_n$ consists of the braids whose strands begin and end at the same position.  Now fix any set of $n$ objects $X_{i(1)},\ldots,X_{i(n)}$.  Then $\PP_n$ acts on $\End(\bigotimes_j X_{i(j)})$ in the obvious way using the braiding operators of the form $(C_{X,Y})^2$ and their conjugates.  One might ask if the image of $\PP_n$ is finite or infinite for all $n$ and all choices of $X_{i(j)}$.  But this is not a more general question: If we define $X=\bigoplus_jX_{i(j)}$ then if $\CC$ has property $F$, $\B_n$ has finite image on $\End(X^{\ot n})$, so that by restricting to $\PP_n$ and to the subspace $\End(\bigotimes_j X_{i(j)})\subset\End(X^{\ot n})$, one sees that the $\PP_n$ image is finite as well.  Obviously the converse is true as well: since $\PP_n$ has finite index, we may take $X_{i(j)}=Y$ for all $j$ and so finiteness of the $\PP_n$ image implies finiteness of the $\B_n$ image.  Similar statements can be made if we replace $\PP_n$ by any finite index subgroup of $B_n$ obtained as a pull-back of a subgroup of $S_n$ via the homomorphism above.  For example, the subgroup of $\B_n$ generated by those elements with the first strand beginning and ending at the same vertical position is the pull-back of the subgroup of $S_n$ that fixes $1$.
\end{remark}

\subsection{Link Invariants}
Associated to any modular category $\CC$ is a $3D$-TQFT, which gives rise to $3$-manifold and link invariants.  In
essence the link invariants are obtained by representing a link $L$ as the closure of a braid $\beta\in\B_n$, and then taking the trace of the image of $\beta$ in one of the representations $\phi_X^n$ of $\B_n$ described above.  More generally, one colors each component of $L$ with objects $X_{i(j)}$ of $\CC$ and represents the colored link as the closure of a braid $\gamma$ where the strands of $\gamma$ must respect the given coloring.  Then one takes the trace of the image of $\gamma$ in the appropriate endomorphism space.  See \cite[Chapter II]{Tur} for full details.  There are two standard choices that will appear below.  We consider the invariants corresponding to coloring all components with either a fixed simple object $X_i$ or the sum of all simple objects.

The link invariants associated to the modular categories mentioned above are as follows, where $q=e^{\pi i/\ell}$:
\begin{enumerate}
\item The link invariant associated to $\CC(\fsl_2,q)$ where we color each component with the object analogous to the irreducible $2$-dimensional representation of $\fsl_2$ is the
Jones polynomial $J_L(q^2)$.
\item More generally, the link invariant associated to $\CC(\fsl_n,q)$
is a one-variable specialization of the (reparameterized) HOMFLYPT polynomial $P^\prime_L(q,n)$.  As above, the invariant $P^\prime_L$ corresponds to coloring each strand with 
the object analogous to the $n$-dimensional representation of $\fsl_n$.  In the setting of Hecke algebras, this corresponds to the $n$-row quotient.
\item Consider the category $\CC(\g,q)$ where $\g$ is of Lie type $B,C$ or $D$, and in the first two cases $\ell$ is even, and let $X$ be the object analogous to the vector representation of $\g$.  Then the invariant associated
to $X$ is a specialization of the (Dubrovnik version) of the Kauffman polynomial $F_L(q^k,q)$, where $k$ depends on the rank of $\g$.
\item Link invariants associated with $\CC(\g,q)$ for $\g$ of other Lie types have not been extensively studied, nor have invariants associated with
objects other than those analogous to the vector representation.  There are two exceptions. Explicit skein relations have been
worked out by G. Kuperberg for Lie type $G_2$.  Also, the invariant associated with
the object analogous to the fundamental spin representation of $\mathfrak{so}_p$ in $\CC(\mathfrak{so}_{p},q)$ with $\ell=2p$, $p$ an odd prime is known to
be related to the homology modulo $p$ of the double cyclic cover $M_L$ of $S^3$ branched over the given link $L$ (see \cite{dBG} and \cite{GoldJones}).
 \item The link invariants associated to the modular categories $\Rep(D^\omega G)$
with $\omega=0$ are described in \cite{FQ}.  Specifically, if we color each component of $L$ with the sum of the simple objects (or with $DG$ itself), one gets (a normalization of) the classical link invariant
$$H_L(G)=|\Hom(\pi_1(S^3\setminus L),G)|.$$ That is, for a fixed link $L$ it counts the homomorphisms from the fundamental group of the
link-complement to the finite group $G$.
\item the TQFTs associated with doubled spherical categories are usually called
Turaev-Viro(-Ocneanu) TQFTs.  The associated link invariants are not well-studied, although some attention has been paid to two of these ``exotic'' examples, see \cite{HRW}.
\end{enumerate}

Later the computational complexity of evaluating these link invariants will be discussed.  Two important complexity classes are $FP$ and $\#P$.  The class of functions that are computable in polynomial time in the length of the input are of complexity $FP$, which is most closely associated with decision problems of complexity $P$.  The class of counting functions of complexity $\#P$ are related to decision problems of complexity $NP$, where instead of asking if there exists a ``yes'' answer one counts the number of ``yes'' answers.  For example, deciding if a given Boolean expression $E$ has an assignment of truth values that satisfy $E$ is $NP$-complete, while counting the number of such assignments is $\#P$-complete.

\section{The Paradigms}
The two paradigms are shown in Figures \ref{fig:dense} and \ref{fig:finite}
respectively.  Each has three boxes representing braid group images, complexity of link invariants, and utility in quantum computation.  Our limited expertise in physics led us to exclude any corresponding speculations from the paradigm, however, see Remark \ref{conclusions} below.

\subsection{Dense Image Paradigm}
\begin{figure}[t0]

\centerline{\includegraphics[width=3.45in]{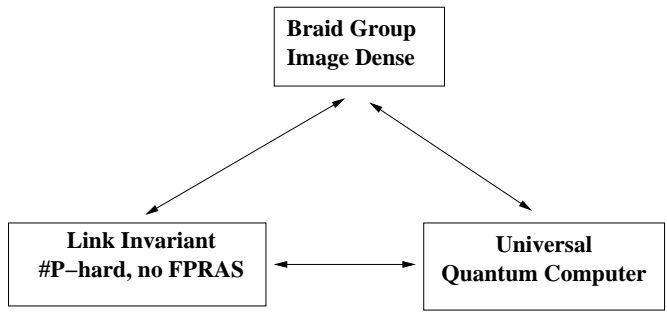}}

\caption{Dense Image Paradigm}

\label{fig:dense}

\end{figure}
In Figure \ref{fig:dense}, ``Braid group image dense'' represents those unitary modular categories which have the density property.

The ``Link invariant'' box requires some explanation.  We say computation of the link invariants are $\#P$-hard because in each known case the exact computation of the invariant can be reduced to a counting problem.  For example, an evaluation of the Jones polynomial of a link $L$ at a root of unity $q^2$ is an integer linear combination of the Galois conjugates of $q$, so that computing each coefficient
may be regarded as a counting problem.  That such an evaluation is \emph{hard} means that if we could find an efficient algorithm for such a problem, we could (in principle) adapt our algorithm to efficiently solve any $\#P$ problem.  However, comparing the quantum computation of a link invariant to classical exact computation is at some level unrealistic for at least two reasons: 1) quantum computation is probabilistic, while classical computation is deterministic and 2) most quantum computations will involve approximate application of some quantum gate (unitary operator), so that the output will be an approximate evaluation as well.  A more relevant question to ask is: does a link invariant $f$ have a \emph{fully polynomial randomized approximation scheme} (FPRAS)?
That is, does there exist an algorithm whose input is a link $L$ with braid index at most $n$ and an error threshold $\ve>0$, whose output is a number $Y$ so that
$$Pr\left(\frac{1}{1+\ve}\le \frac{Y}{f(L)}\le 1+\ve\right)>3/4$$ that runs in polynomial time
in $n$ and $1/\ve$?  Of course by running such an algorithm multiple times, one may improve the certainty that the approximation of $f(L)$ is correct within an $\ve$ factor of $f(L)$.  The associated decision problem complexity class is $RP$ (\emph{radomized polynomial time}).  It is widely believed that $RP\neq NP$, and the non-existence of an $FPRAS$ for a given problem is usually proved under this assumption.

\subsection{Finite Image Paradigm}

\begin{figure}[t0]

\centerline{\includegraphics[width=3.45in]{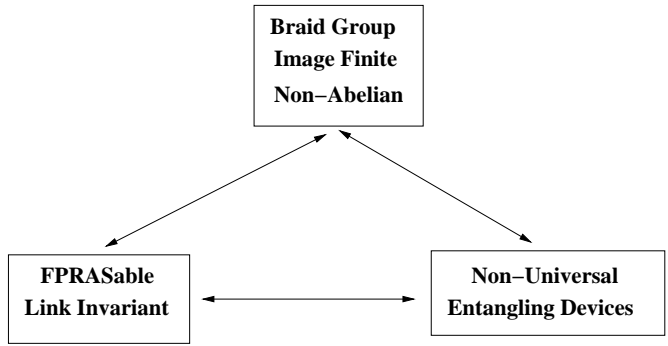}}

\caption{Finite Image Paradigm}

\label{fig:finite}

\end{figure}

Most of the relationships in the Finite Image paradigm (Figure \ref{fig:finite}) can be understood from the remarks on the Dense Image paradigm above.  Notice that we have excluded the finite abelian braid group images from the description.  This is because the cases where the images of the braid group are finite abelian are mathematically trivial, corresponding to \emph{abelian anyons}.  Firstly, the link invariant will essentially count components or at
best linking numbers, which can be done classically in polynomial time.  Secondly, the representations of the braid group in these cases are all $1$-dimensional.  Because of this, there is no ground state degeneracy and hence any device based upon such systems would not even be capable of efficiently storing information, i.e. they would be non-entangling.  It is interesting to note that, to date, the only topological phases that have been convincingly shown to exist are abelian anyons (see \cite{dSetal}).

Non-universal quantum devices that can at least produce entangled qubits could potentially be used to store quantum information and even be useful in quantum error correction (see e.g. \cite{YRWGW}).  A well-known example is the Bell basis change matrix which is related to the Jones polynomial at $t=i$.
\subsection{Evidence}
Analyses of the braid group images and the computational complexity of the link invariant evaluations associated to many of the modular categories described above have been carried out.  We discuss each in turn, recording the precise evidence
for the two paradigms in Table \ref{evidence} where speculations are in bold type.  For notational convenience, set $q=e^{\pi i/\ell}$, and denote by $c$ the number of components of a link $L$.  Let $d_k$ be the dimension of the homology space modulo $k$ of the double cyclic cover of $S^3$ branched over $L$.  $K(L)$ is a classical invariant that only depends on the linking matrix of $L$.

\subsubsection{\textbf{Jones polynomial}}
For $\CC(\fsl_2,q)$ with $q=e^{\pi i/\ell}$, $3\leq\ell$ and $X_1$ the object corresponding to the fundamental $2$-dimensional representation of $\fsl_2$ the algebra $\End(X_1^{\ot n})$ is isomorphic to the Temperley-Lieb algebra $TL_n(q^2)$.  Jones determined precisely when the braid group images are finite in \cite{jones86}, and in all other cases it is shown in \cite[Theorem 0.1]{FLW2} that the braid group images are dense.  The (exact) computational complexity of the corresponding link invariant the Jones polynomial $J_L(q^2)$ was worked out in \cite{JVW}, where it is shown that, except for $\ell\in\{1,2,3,4,6\}$ the complexity class is $\#P$-hard.  This was accomplished by using a result of Thistlethwaite that evaluating the Jones polynomial at $t=q^2$ for $L$ an alternating link is essentially equivalent to computing the Tutte polynomial of an associated plane graph $G(L)$ at $(-t,-1/t)$, which is shown to be $\#P$-hard except at the special points described above.  The Jones polynomial at these special points degenerates to a ``classical'' link-invariant that is computable in polynomial time.  We conjecture the following:

\begin{conjecture}
 There is no FPRAS for evaluating $J_L(q^2)$ except at the special points described above, provided $RP\neq NP$.
\end{conjecture}
This conjecture is partially motivated by the belief that quantum computers are strictly more powerful than classical computers.  If this conjecture were false, there would be an $FPRAS$ for a $BQP$-complete problem.  A second, less philosophical, piece of evidence is found in \cite{GJ}, where it is shown that, away from the positive quadrant in the rational $xy$-plane and a few exceptional curves, no $FPRAS$ exists for evaluating the Tutte polynomial at $(x,y)$.  This result does not apply to complex pairs $(x,y)$ and so does not give any information for the Jones polynomial at roots of unity, but is nonetheless compelling evidence for our conjecture.

\subsubsection{\textbf{HOMFLYPT polynomial}}
Generalizations of the results above to the categories $\CC(\fsl_n,q)$ with $q=e^{\pi i/\ell}$ and the corresponding specializations of the the HOMFLYPT polynomial are found in \cite{Welsh} (due to Vertigan), \cite{FLW2}, and \cite{GoldJones}.  The role of the Temperley-Lieb algebra is taken by specializations of the two-parameter Hecke-algebra (see \cite{jones87}), and the results are of the same format with one exception: for $n\geq 3$ one may have infinite braid group images that are not dense, see \cite[Theorem 4.1]{FLW2}.
See also \cite{MOO} for a related invariant obtained by summing over all simple objects.

Since the Jones polynomial can be obtained as a specialization of the  HOMFLYPT polynomial, FPRASability of the HOMFLYPT polynomial would imply the same for the Jones polynomial.

\subsubsection{\textbf{Kauffman polynomial}}
The computational complexity of evaluating of the Kauffman polynomial has been worked out by Vertigan, see \cite{Welsh}.  The relevant modular categories are obtained from the categories of the form $\CC(\g,q)$ with $\g\in\{\mathfrak{so}_N,\mathfrak{sp}_{2N}\}$.  In these cases the algebras $\End(X^{\ot n})$ are related to specializations of the form $r=q^k$ of $BMW$-algebras $C_n(r,q)$ (see \cite[Prop. 2.1]{wenzlsurvey}), where as usual $X$ is the quantum analogue of the vector representation.  The braid group images are worked out in all non-trivial cases except $r=\pm i$ in
\cite{J2},
\cite{LRW}, \cite{LR} and \cite{Jthesis}.  In general the images are either finite or dense, although exceptions are found in \cite{LR}, and are expected for $r=\pm i$.

Again, as the Jones polynomial can be obtained as a specialization of the Kauffman polynomial, FPRASability of the Kauffman polynomial would imply the same for the Jones polynomial.

\subsubsection{$\mathbf{d_n=\dim H_1(M_L,\Z_n)}$} The categories $\CC(\mathfrak{so}_{2n+1},e^{\pi i/\ell})$ with $\ell=2(2n+1)$ may be regarded as the extension of the series of modular categories whose first two terms are $\CC(\fsl_2,e^{\pi i/6})$ and $\CC(\mathfrak{sp}_4,e^{\pi i/10})$.  At least for $2n+1=p\geq 7$ prime, the corresponding link invariants  are $\pm(\sqrt{p})^{d_p}$ where $d_p=\dim H_1(M_L,\Z_p)$ with $M_L$ the double cyclic cover of $S^3$ branched over $L$, see \cite{GoldJones} and \cite{dBG}.  Polynomial algorithms exist for computing the dimension of these homology spaces, and the braid group images are shown to be finite (symplectic) groups in \cite{GoldJones}.  It seems reasonable that this should hold for arbitrary $2n+1$
as well.

\subsubsection{$\mathbf{H_L(G)=|\Hom(\pi_1(S^3\setminus L),G)|}$}
While the fact that the invariant corresponding to the modular category $\Rep(DG)$ for $G$ a finite group is the classical invariant $H_L(G)$ has been known for some time, the computational complexity has not been studied to our knowledge.  Moreover, the fact that $\Rep(DG)$ has property $F$ was shown only recently \cite{ERW}.

Recent results suggest the following:
\begin{conjecture} Let $G$ be a finite group and $L$ a link.
 \begin{enumerate}
  \item[(a)] There exists an $FPRAS$ for computing $H_L(G)$ for any group $G$.
\item[(b)] Suppose $G$ is solvable.  Then there is a polynomial algorithm for exact computation of $H_L(G)$.
 \end{enumerate}
\end{conjecture}
We support this conjecture with the following list of facts:
\begin{enumerate}
 \item Clearly if $G$ is an abelian group and $L$ has $k$ components, then $H_L(G)=|\Hom(H_1(S^3\setminus L),G)|=|\Hom(Z^k,G)|=|G|^k$.
\item It is shown in \cite{Eis} that if $G$ is nilpotent and $L$ is a knot
then $H_L(G)=|G|$ is constant.  So at least for knots, $H_L(G)$ is polynomial time computable for $G$ nilpotent.
\item In \cite{MS} an algorithm for computing the number of homomorphisms from
a given finitely presented group $\Gamma$ to a finite solvable group is given.  It is not
clear if this algorithm finishes in polynomial time (in, say, the number of generators of $\Gamma$), but it certainly supports the case for (b).  Moreover, in preliminary computations (worked out with S. Witherspoon) for $G$ a generalized dihedral group we found that the corresponding braid group representation is equivalent to a finite field evaluation of the Burau representation.  This is significant, as the Burau representation supports the
Alexander polynomial, which is known to be polynomial-time computable.
\item Even in the non-solvable case, an algorithm exists: $\pi_1(S^3\setminus L)$
has presentation $\lan x_1,\ldots,x_n: R_1,\ldots,R_m\ra$ with $n$ and $m$ are
bounded by $N+M$ where $N$ is the number of strands in some projection of $L$ and $M$ is the number of crossings.  One checks all $|G|^n$ $n$-tuples against the $m$ relations to find homomorphisms.  One could improve this algorithm slightly by applying automorphisms of $G$, but the algorithm would still be exponential in $n$.  Perhaps a randomization of this algorithm where one samples a moderately-sized subset of the $n$-tuples of elements of $G$ and then approximates $H_L(G)$ by proportion would provide an $FPRAS$.  Whether this could be done efficiently and accurately would require some analysis.  We should mention that it is widely believed that an $FPRAS$ exists for computing $H_L(G)$ (\cite{Kup}).

\end{enumerate}

\begin{table}\label{evidence}
\begin{tabular}{*{2}{|p{1.9cm}}|p{3.2cm}|p{2.2cm}|p{1.75cm}|}
\hline
\rule[-2mm]{0mm}{6mm}UMC & Restrictions & Invariant & Complexity & $\B_n$ Image \\
\hline\hline
\rule[-2mm]{0mm}{6mm}$\CC(\fsl_2,q)$ & $5\leq\ell\neq 6$ & $V_L(q^2)$ & \raggedright{$\#P$-hard \\  \textbf{no FPRAS?}} & dense  \\
\hline
\rule[-2mm]{0mm}{6mm}\raggedright{$\CC(\mathfrak{sl}_n,q)$,  $3\leq n$} & \raggedright{$n+2\leq\ell$,\\ $\ell\neq 6$}
& $P^\prime_L(q,n)$  & \raggedright{$\#P$-hard \\  \textbf{no FPRAS?}}& infinite \quad
not dense  \\
\hline
\rule[-2mm]{0mm}{6mm}\raggedright{$\CC(\mathfrak{so}_{2n+1},q)$, $2\leq n$} & \raggedright{$\ell$ even,\\ $2n+2\leq\ell$,\\ $\ell\not=4n$} & $F_L(q^{2n},q)$ & \raggedright{$\#P$-hard \\  \textbf{no FPRAS?}}& dense  \\
\hline
\rule[-2mm]{0mm}{7mm}\raggedright{$\CC(\mathfrak{sp}_{2n},q)$, $2\leq n$} & \raggedright{$\ell$ even,\\ $2n+6\leq\ell$,\\ $\ell\not=4n+2$} & $F_L(q^{-2n-1},q)$ & \raggedright{$\#P$-hard \\  \textbf{no FPRAS?}} & dense  \\
\hline
\rule[-2mm]{0mm}{7mm}\raggedright{$\CC(\mathfrak{so}_{2n},q)$, $3\leq n$}& \raggedright{$2n+2\leq\ell$,\\ $\ell\not=4n-2$} & $F_L(q^{2n-1},q)$ & \raggedright{$\#P$-hard \\  \textbf{no FPRAS?}} & dense  \\
 \hline \rule[-2mm]{0mm}{7mm}$\CC(\mathfrak{so}_{4},q)$& $7\leq\ell$ & $(-1)^{c-1}[V_L(-q^{-2})]^2$ & \raggedright{$\#P$-hard \\  \textbf{no FPRAS?}} & infinite \quad
not dense \\
 \hline
\rule[-2mm]{0mm}{7mm}$\CC(\fsl_2,q)$ & $\ell=3$ & $(-1)^{c-1}$  & $FP$  & finite\quad abelian  \\
 \hline
\rule[-2mm]{0mm}{7mm}$\CC(\fsl_2,q)$ & $\ell=4$ & $(-\sqrt{2})^{c-1}(-1)^{\text{Arf}(L)}$\quad or $0$  & $FP$  & finite \\
\hline
\rule[-2mm]{0mm}{7mm}$\CC(\fsl_n,q)$ & $\ell=6$ & $\pm(i)^{c-1}(i\sqrt{3})^{d_3}$  & $FP$  & finite  \\
 \hline
\rule[-2mm]{0mm}{7mm}$\CC(\mathfrak{sp}_4,q)$ & $\ell=10$ & $\pm(\sqrt{5})^{d_5}$  & $FP$  & finite  \\
\hline

\rule[-2mm]{0mm}{7mm}$\CC(\mathfrak{sl}_n,q)$ & $\ell=n+1$  & $e^{\pi i K(L)/n}$ & $FP$ & finite \quad abelian  \\
\hline
\rule[-2mm]{0mm}{7mm}\raggedright{$\CC(\mathfrak{so}_{p},q)$,\quad $3\leq p$ prime}\\ & $X$ spin rep., $\ell=2p$ & $\pm(\sqrt{p})^{d_{p}}$ & $FP$ & finite  \\
\hline
\rule[-2mm]{0mm}{7mm}$\Rep(DG)$ & $G$ finite & $H_L(G)$  & \textbf{FPRAS?} & finite  \\

\hline
\end{tabular}
\end{table}

\begin{remark}\label{conclusions}
 \begin{enumerate}
 \item
 We speculate that an appropriate physical aspect of the paradigm would be as follows:  When the braid group image is dense, then it is unlikely that there is
an efficient way to approximately simulate the corresponding physical system.
Our expertise in the subject is not sufficient to say anything authoritative, but it seems reasonable that an efficient approximate simulation of the physical system could be used to construct an FPRAS for the link invariant.
When the braid group image is finite, we might expect that efficient numerical methods (such as quantum Monte Carlo) exist for (approximately) simulating the corresponding quantum mechanical systems (see e.g. \cite{JOVVC}).
\item There is a related conjecture characterizing UMCs with property $F$ by the categorical dimensions of their
simple objects.  This is beyond our current scope, but details will appear in \cite{propF}.
\end{enumerate}
\end{remark}

\pagebreak


\begin{thebibliography}{99}


\bibitem{BK} B.\ Bakalov; A.\ Kirillov, Jr., {\em Lectures on
Tensor Categories and Modular Functors}, University Lecture
Series, vol.\ {\bf 21},  Amer.\ Math.\ Soc., 2001.
\bibitem{BarWest} J.\ W.\ Barrett; B.\ W.\ Westbury, Spherical categories.  \textit{Adv. Math.}  \textbf{143}  (1999),  no. 2, 357--375.
\bibitem{dSetal} S.\ Das Sarma; M.\ Freedman; C.\ Nayak; S.\ H.\ Simon,; A.\ Stern,
Non-Abelian Anyons and Topological Quantum Computation.
arXiv:0707.1889.

\bibitem{dBG} J.\ de Boer; J.\ Goeree, Markov traces and ${\rm II}\sb 1$ factors in conformal field theory.  \textit{Comm. Math. Phys.}  \textbf{139}  (1991),  no. 2, 267--304.
\bibitem{DW} R.\ Dijkgraaf; E.\ Witten, \emph{Topological gauge theories and group cohomology.}
 Comm. Math. Phys. 129 (1990), no. 2, 393--429.
\bibitem{Eis} M.\ Eisermann, The number of knot group representations is not a Vassiliev invariant.  \textit{Proc. Amer. Math. Soc. } \textbf{128}  (2000),  no. 5, 1555--1561.
\bibitem{ERW} P.\ Etingof; E.\ C.\ Rowell; S.\ J.\ Witherspoon,
Braid group representations from quantum doubles of finite groups. \emph{Pacific J. Math.} \textbf{234}, no. 1 (2008), 33-41.
\bibitem{Jthesis} J.\ Franko, Braid group representations via the Yang Baxter Equation, Ph.D. thesis, Indiana University, 2007.
\bibitem{FRW}  J.\ Franko; E.\ C.\ Rowell; Z.\ Wang, Extraspecial 2-groups and images of braid group
representations, \emph{J. Knot Theory Ramifications} \textbf{15} (2006) no. 4, 1-15.

\bibitem{FQ} D.\ Freed; F.\ Quinn, Chern-Simons theory with finite gauge group.
\textit{ Comm. Math. Phys.} \textbf{156} (1993), no. 3, 435--472.

\bibitem{FPNAS}Freedman, Michael H. P/NP, and the quantum field computer.  Proc. Natl. Acad. Sci. USA  95  (1998),  no. 1, 98--101 (electronic).
\bibitem{FKLW}M.\ Freedman; A.\ Kitaev; M.\ Larsen; Z.\ Wang,
Topological quantum computation.
\textit{Bull. Amer. Math. Soc.
(N.S.)} 40 (2003), no. 1, 31--38.
\bibitem{FLW2} M.\ H.\ Freedman; M.\ J.\ Larsen; and Z.\ Wang,
The two-eigenvalue problem and density of Jones representation of
braid groups. \textit{Comm. Math. Phys.} \textbf{228} (2002),
177-199, arXiv: math.GT/0103200.
\bibitem{GJ} L.\ A.\ Goldberg; M.\ Jerrum, Inapproximability of the Tutte polynomial, in \emph{STOC '07: Proceedings of the thirty-ninth annual ACM symposium on Theory of computing}, 459--468, ACM, New York 2007.
\bibitem{GoldJones} D.\ M.\ Goldschmidt; V.\ F.\ R.\ Jones, Metaplectic link invariants.  \textit{Geom. Dedicata } \textbf{31}  (1989),  no. 2, 165--191.
\bibitem{HRW} S.-M.\ Hong; E.\ C.\ Rowell; Z.\ Wang, On exotic modular tensor categories. arXiv:0710.5761.
\bibitem{Izumi}M.\ Izumi, The structure of sectors associated with
Longo-Rehren inclusions. II. Examples. \textit{Rev. Math. Phys.} \textbf{13}
(2001), no. 5, 603--674.

\bibitem{JVW} F.\ Jaeger; D.\ L.\ Vertigan; D.\ J.\ A.\ Welsh, On the computational complexity of the Jones and Tutte polynomials.  \textit{Math. Proc. Cambridge Philos. Soc. } \textbf{108}  (1990),  no. 1, 35--53.
\bibitem{jones86} V. F. R. Jones, Braid groups, Hecke algebras and type ${\rm II}\sb 1$ factors,
\textit{Geometric methods in
 operator algebras (Kyoto, 1983)}, 242--273, Pitman Res. Notes Math. Ser., 123, Longman Sci. Tech.,
 Harlow, 1986.
\bibitem{jones87} V.\ F.\ R.\ Jones, Hecke algebra representations of braid groups and link polynomials.  \textit{Ann. of Math. (2)}  \textbf{126}  (1987),  no. 2, 335--388.
\bibitem{J2} V. F. R. Jones, On a certain value of the Kauffman polynomial.
\textit{Comm. Math. Phys.} \textbf{125} (1989), no. 3, 459--467.

\bibitem{JOVVC} J.\ Jordan; R.\ Orus; G.\ Vidal; F.\ Verstraete; J.\ I.\ Cirac, Classical simulation of infinite-size quantum lattice systems in two spatial dimensions
arXiv: cond-mat/0703788.

\bibitem{kauf} L.\ Kauffman, An invariant of regular isotopy. \emph{Trans. Amer. Math. Soc.} \textbf{318} (1990), no. 2, 417--471.
\bibitem{KL} L.\ Kauffmann, S.\ Lomonaco Jr., Braiding
operators are universal quantum gates. \textit{New J. Phys.} \textbf{6} (2004), 134.1-134.40 (electronic).
\bibitem{Kit} A.\ Kitaev, {Fault-tolerant quantum computation by anyons}.
\textit{Ann. Physics } \textbf{303} (2003), no. 1, 2--30.
\bibitem{Kup} G.\ Kuperberg, private communication.

\bibitem{LR}  M.\ J.\ Larsen; E.\ C.\ Rowell, An algebra-level version of a link-polynomial identity of Lickorish, to appear in \emph{Math.
Proc. Cambridge Philos. Soc.}
\bibitem{LRW} M.\ J.\ Larsen; E.\ C.\ Rowell; Z. Wang, The $N$-eigenvalue problem and two applications.  \textit{Int. Math. Res. Not.} \textbf{2005} (2005), no. 64, 3987--4018.

\bibitem{Li} W.\ B.\ R.\ Lickorish, Some link-polynomial relations.  \textit{Math. Proc. Cambridge Philos. Soc.} \textbf{105}  (1989),  no. 1, 103--107.

\bibitem{Ll} S.\ Lloyd, Quantum computation with abelian anyons.  Quantum Inf. Process.  1  (2002),  no. 1-2, 13--18.
\bibitem{MS} D.\ Matei; A.\ I.\ Suciu, Counting homomorphisms onto finite solvable groups.  \textit{J. Algebra}  \textbf{286}  (2005),  no. 1, 161--186.
\bibitem{Muger2} M.\ M\"uger, \emph{From subfactor to categories and topology, II}
J. Pure Appl. Algebra 180 (2003), no. 1-2, 159--219.

\bibitem{MOO} H.\ Murakami; T.\ Ohtsuki; M.\ Okada, Invariants of three-manifolds derived from linking matrices of framed links. \textit{ Osaka J. Math.}  \textbf{29}  (1992),  no. 3, 545--572.
\bibitem{Rowell1} E. C. Rowell: On a family of non-unitarizable ribbon categories:  \textit{Math Z.}
\textbf{250} no. 4 (2005) 745--774.

\bibitem{Rsurvey} E.\ C.\ Rowell From quantum groups to unitary modular tensor categories
in \textit{Contemp.\ Math.\ } \textbf{413} (2006), 215--230.
\bibitem{RJPAA} E.\ C.\ Rowell, \emph{Unitarizablity of premodular
categories} preprint, arXiv: 0710.1621, to appear in J. Pure Appl. Algebra.

\bibitem{propF} E.\ C.\ Rowell, \emph{A finiteness property for braided
tensor categories} in preparation.


\bibitem{Tur}V.\ Turaev, Quantum Invariants of Knots and 3-Manifolds, De Gruyter Studies in
Mathematics, Walter de Gruyter 1994.

\bibitem{Welsh} D.\ Welsh,
  \emph{Complexity: Knots, Colourings and Counting}, LMS Lecture
  Notes Series 186, Cambridge University Press, Cambridge, 1993.
\bibitem{wenzlsurvey} H.\ Wenzl, Tensor categories and braid representations, in  \textit{Quantum groups and Lie theory (Durham, 1999)},  216--234, London Math. Soc. Lecture Note Ser., 290, Cambridge Univ. Press, Cambridge, 2001.
\bibitem{Wenzlcstar}H.\ Wenzl, $C\sp *$ tensor categories from quantum groups.
\textit{ J. Amer. Math. Soc.} \textbf{11} (1998), no. 2, 261--282.
\bibitem{YRWGW} Y.\ Zhang; E.\ C.\ Rowell; Y.-S.\ Wu; Z.\ Wang; M.-L.\ Ge, From extraspecial two-groups to GHZ states. arXiv:0706.1761.
\end{thebibliography}
\end{document}